# A Statistical Approach to Broken Stick Problems

Rahul Mukerjee
Indian Institute of Management Calcutta
Joka, Diamond Harbour Road
Kolkata 700 104, India
e-mail: rmuk0902@gmail.com

*Abstract*: Let a stick be broken at random at $n-1$ points to form $n$ pieces. We consider three problems on forming $k$-gons with $k$ out of these $n$ pieces, $3 \leq k \leq n$, and show how a statistical approach, through a linear transformation of variables, yields simple solutions that also allow fast computation.

**1. Introduction**

Consider a stick broken at random at $n-1$ points to form $n$ pieces. We explore three problems on the formation of $k$-gons with $k$ out of these $n$ pieces, where $3 \leq k \leq n$:

I. What is the probability that no choice of $k$ pieces out of $n$ forms a $k$-gon?

II. What is the probability that every choice of $k$ pieces out of $n$ forms a $k$-gon?

III. If $k$ pieces are chosen at random out of the $n$ pieces, all such choices being equally probable, then what is the probability that the $k$ chosen pieces form a $k$-gon?

Although the broken stick and related problems have a long history (see e.g., Clifford, 1866; D'Andrea and Gómez, 2006; Goodman, 2008), it is only very recently that problems I and II received attention in full generality. Verreault (2022a) studied problem II and found a unified analytical solution that covers earlier results pertaining to special cases of $k$ and $n$. This was done via a link with order statistics from an exponential distribution. He also mentioned I as an open problem and remarked that it remains intractable even if one attempts to use the link with exponential order statistics. Thereafter, Verreault (2022b) obtained a solution to problem I using deep combinatorial techniques which first consider integer valued breaking points and associated Diophantine inequalities, then employs MacMahon's (1915-16) partition analysis, and finally passes from the discrete to the continuous case. His ultimate interest was in the continuous case and his main theorem, indeed, dwells on this case, giving a solution to problem I via generalized Fibonacci numbers. We refer to Verreault (2022a, b) for further references; see these papers the references there for an indication of how the broken stick model can be useful in diverse areas ranging from biology (MacArthur, 1957) to finance (Tashman and Frey, 2009).

In this paper, one of our main objectives is to report an alternative solution to problem I through a direct statistical approach. Specifically, we show that the issue left open in Verreault (2022a) regarding this problem can be settled quickly by exploiting the link with exponential order statistics to a fuller extent by means of a linear transformation of these statistics. Our solution, presented in Section 2, is much simpler and shorter than the one in Verreault (2022b). We hope that this will be appealing to a statistical audience because it shows the power of relatively straightforward statistical considerations over more demanding combinatorial arguments in so far as a resolution of problem I is concerned. For academic interest, we also satisfy ourselves in the appendix that, with appropriately defined notation, our solution perfectly matches the one in Verreault (2022b).

The idea of linear transformation employed for problem I is pursued in Section 3, where we show that it yields an appreciably shorter solution also to problem II, compared to Verreault (2022a). Finally, problem III is investigated in Section 4, leading to a surprise answer which depends on $k$ alone, and not $n$.

For ease in reference, before concluding the introduction, we briefly describe the aforesaid link between the broken stick problem and exponential order statistics. Without loss of generality, let the stick have length 1, and let $\Delta_1,...,\Delta_n$ be the lengths of the $n$ pieces when it is broken at random at $n$ – 1 points. Write $\Delta_{(1)} < ... < \Delta_{(n)}$ for the order statistics corresponding to $\Delta_1,...,\Delta_n$. Also, let $Y_1 < \ldots < Y_n$ be the order statistics based on a random sample of size $n$ from the exponential distribution with mean 1, and $W = Y_1 + \ldots + Y_n$. From earlier findings in Rényi (1953) on order statistics (see also Pyke, 1965), Verreault (2022a) deduced the following proposition.

**Proposition 1**. *The joint distribution of* $\Delta_{(1)},...,\Delta_{(n)}$ *is the same as that of* $Y_1/W,...,Y_n/W$.

## 2. Problem I: Probability of never forming a *k*-gon

Let $e_1,\ldots, e_n$ be the unit vectors of order $n$x1, and define vectors $b_1,\ldots, b_n$ recursively as

$$b_1 = e_1,\ b_r = e_r + b_{r-1}\ (r = 2,\ldots, k-1), \quad b_r = e_r + \Sigma_{u=1}^{k-1} b_{r-u}\ (r = k,\ldots, n). \tag{1}$$

Let $\Sigma_{r=1}^{n} b_r = (\beta_1,...,\beta_n)'$, where the prime stands for transpose. Then the following result, giving a simple solution to Problem I through a direct statistical approach, holds.

**Theorem 1**. *The probability that no choice of k pieces out of n forms a k-gon is given by*

$$P(k,n) = n!/(\Pi_{r=1}^{n}\beta_r).$$

*Proof*. By Proposition 1,

$$P(k,n) = P(Y_r > \Sigma_{u=1}^{k-1} Y_{r-u},\ r = k,\ldots, n);$$

cf. Verreault (2022a, Section 3). Thus,

$$P(k,n) = n!\int_S \exp(-\Sigma_{r=1}^{n} y_r)dy, \tag{2}$$

where $dy \equiv dy_1\ldots dy_n$ and

$$S = \{(y_1,\ldots, y_n): 0 < y_1 < \ldots < y_{k-1} \text{ and } y_r > \Sigma_{u=1}^{k-1} y_{r-u},\ r = k,\ldots, n\}. \tag{3}$$

Make a linear variable transformation, with Jacobian unity, as given by

$$x_1 = y_1,\ x_r = y_r - y_{r-1}\ (r = 2,\ldots, k-1),\quad x_r = y_r - \Sigma_{u=1}^{k-1} y_{r-u}\ (r = k,\ldots, n). \tag{4}$$

by (3) and (4), $(y_1,\ldots, y_n) \in S$ if and only if $x_1,\ldots, x_n > 0$. Moreover, writing $x = (x_1,\ldots, x_n)'$, by (1) and (4), $y_r = b_r' x$, $(r = 1,\ldots, n)$, so that $\Sigma_{r=1}^n y_r = \Sigma_{r=1}^n b_r' x = \Sigma_{r=1}^n \beta_r x_r$. The result is now immediate from (2). □

Theorem 1, aided by the recursion relation in (1), yields $P(k,n)$ almost instantaneously. This happens even for relatively large $n$ and $k$ such as $(n, k) = (200, 196)$, when it is advisable to apply the theorem as $P(k,n) = \Pi_{r=1}^n \{(n-r+1)/\beta_r\}$ to avoid division of a big number by another. Of course, $P(k,n)$ is typically close to zero for large $n$, a point that one expects intuitively and is foreshadowed by the fact that $P(n,n) = n/2^{n-1}$; see Verreault (2022a).

While Theorem 1 alone suffices for the practical purpose of computing $P(k,n)$ efficiently, it is of considerable theoretical interest to reconcile it with the corresponding result in Verreault (2022b). This is done in the appendix based on a backward recursion relation among $\beta_1,\ldots,\beta_n$, as given by

$$\beta_n = 1,\qquad \beta_j = 1 + \Sigma_{r=j+1}^{j+k-1} \beta_r\quad (j = k-2,\ldots, n-1),$$

$$\beta_j = 1 + \beta_{j+1} + \Sigma_{r=k}^{j+k-1} \beta_r\qquad (j = 1,\ldots, k-3), \tag{5}$$

where the last set of equations arises only when $k \geq 4$, and we take

$$\beta_{n+1} = \ldots = \beta_{n+k-2} = 0. \tag{6}$$

The initial conditions (6) make the sums in (5) well-defined even when $j + k - 1$ exceeds $n$. The above recursion relation is proved in the appendix where we also note certain simplifying features related to Verreault's (2022b) result. Together, these quickly lead to the desired reconciliation.

We remark that the recursion relation in (5) and (6) does not entail any computational gain over the direct calculation $(\beta_1,\ldots,\beta_n)' = \Sigma_{r=1}^n b_r$ using (1), because as noted earlier, the latter is almost instantaneous. The only purpose of (5) and (6) is matching Theorem 1 with Verreault's (2022b) result. This reconciliation also shows that even if one is interested in reaching the form in Verreault (2022b) for $P(k,n)$, the statistical approach via the transformation (4) provides a more transparent route, completely bypassing elaborate combinatorial arguments.

## 3. Problem II: Probability of always forming a *k*-gon

In Theorem 2 below, we find a solution to this problem, in terms of the complementary probability, again using the idea of transformation of variables. In what follows,

$$c_{rj} = (j+2)(k-1-r) + n - k + 2\qquad (r = 1,\ldots, k-1;\ j = 0,\ldots, n-k). \tag{7}$$

**Theorem 2**. *The probability that a k-gon cannot be formed for at least one choice of k pieces out of n is given by*

$$\overline{Q}(k,n) = \frac{n!}{(n-k+2)!}\Sigma_{j=0}^{n-k}(-1)^j\binom{n-k+1}{j+1}(\Pi_{r=1}^{k-2}c_{rj})^{-1}.$$

*Proof*. By Proposition 1, $\overline{Q}(k,n) = P(Y_n > \Sigma_{r=1}^{k-1}Y_r)$ ; cf. Verreault (2022a, Section 2). If one recalls the well-known form of the joint density of $Y_1, \ldots, Y_{k-1}$ and $Y_n$ , then one obtains

$$\overline{Q}(k,n) = \frac{n!}{(n-k)!}\int_T \exp\{-(\Sigma_{r=1}^{k-1}y_r + y_n)\}\{\exp(-y_{k-1}) - \exp(-y_n)\}^{n-k} dy$$

$$= \frac{n!}{(n-k)!}\Sigma_{j=0}^{n-k}(-1)^j\binom{n-k}{j}I_j, \qquad (8)$$

using binomial expansion, where $dy \equiv dy_1\ldots dy_{k-1}dy_n$, and

$$I_j = \int_T \exp[-\{\Sigma_{r=1}^{k-1}y_r + (n-k-j)y_{k-1} + (j+1)y_n\}]dy \quad (j = 0,\ldots, n-k), \qquad (9)$$

$$T = \{(y_1, \ldots, y_{k-1}, y_n): 0 < y_1 < \ldots < y_{k-1} \text{ and } y_n > \Sigma_{r=1}^{k-1}y_r\}.$$

Make a linear variable transformation, with Jacobian unity, as given by

$$x_1 = y_1,\ x_r = y_r - y_{r-1}\ (r = 2,\ldots, k-1), \quad x_n = y_n - \Sigma_{r=1}^{k-1}y_r, \qquad (10)$$

which is equivalent to what (4) would become for $n = k$. Then

$$y_r = x_1 + \ldots + x_r\ (r = 1,\ldots, k-1), \quad \Sigma_{r=1}^{k-1}y_r = \Sigma_{r=1}^{k-1}(k-r)x_r$$

and hence

$$\Sigma_{r=1}^{k-1}y_r + (n-k-j)y_{k-1} + (j+1)y_n$$

$$= \Sigma_{r=1}^{k-1}(k-r)x_r + (n-k-j)\Sigma_{r=1}^{k-1}x_r + (j+1)\{\Sigma_{r=1}^{k-1}(k-r)x_r + x_n\} = \Sigma_{r=1}^{k-1}c_{rj}x_r + (j+1)x_n,$$

recalling (7). Moreover, by (10), $(y_1, \ldots, y_{k-1}, y_n) \in T$ if and only if $x_1,\ldots, x_{k-1}, x_n > 0$. Thus, from (9), we obtain $I_j = \{(\Pi_{r=1}^{k-2}c_{rj})(n+k-2)(j+1)\}^{-1}$, as $c_{rj} = n + k - 2$ when $r = k - 1$. The result now follows from (8) after a little simplification. □

The special case of $\overline{Q}(n,n)$, which is the same as $P(n,n)$ in the last section, is of some interest. If $k = n$, then by (7), $c_{r0} = 2(n-r)$ $(r = 1,\ldots, n-2)$ and Theorem 2 reduces to $\overline{Q}(n,n) = n/2^{n-1}$. Indeed, even without invoking Theorem 2, this follows almost effortlessly from the above variable transformation approach. To see this, observe that $\overline{Q}(n,n) = n!\int_T \exp(-\Sigma_{r=1}^{n}y_r)dy$, where $dy \equiv dy_1\ldots dy_n$ and $T = \{(y_1,\ldots, y_n): 0 < y_1 <\ldots< y_{n-1} \text{ and } y_n > y_1 +\ldots+ y_{n-1}\}$. Now, with

$$x_1 = y_1,\ x_r = y_r - y_{r-1}\ (r = 2,\ldots, n-1), \quad x_n = y_n - \Sigma_{r=1}^{n-1}y_r,$$

one gets (i) $\Sigma_{r=1}^{n}y_r = 2\Sigma_{r=1}^{n-1}(n-r)x_r + x_n$ and (ii) $(y_1,\ldots, y_n) \in T$ if and only if $x_1,\ldots, x_n > 0$. So, we obtain $\overline{Q}(n,n) = n!/\{2^{n-1}\Pi_{r=1}^{n-1}(n-r)\} = n/2^{n-1}$, and these arguments appear to be even more elementary than those in Verreault (2022a, Section 2).

For general $k$ and $n$, the proof in Theorem 2 for Problem II is considerably shorter than its counterpart in Verreault (2022a). We now observe that Theorem 2 is in agreement with Verreault's (2022a) main result, according to which the probability that all $k$ pieces form $k$-gons is

$$Q(k,n) = \frac{n!}{(n-k+2)!(n-k+2)} \Sigma_{j=1}^{n-k+2} \{(-1)^{j+1} / j^{k-3}\}\binom{n-k+2}{j}\{(\Pi_{r=0}^{k-3}\lambda_{rj})\}^{-1},$$

where $\lambda_{rj} = r+1+\{(n-k+2)/j\}$. It suffices to show that $\overline{Q}(k,n)$ in Theorem 2 equals $1-Q(k,n)$.

Note that $\Pi_{r=0}^{k-3}\lambda_{r1} = n!/(n-k+2)!$, while for $j \ge 2$,

$$\Pi_{r=0}^{k-3}\lambda_{rj} = (1/j^{k-2})\Pi_{r=0}^{k-3}\{j(r+1)+n-k+2\}$$

$$= (1/j^{k-2})\Pi_{r=1}^{k-2}\{j(k-1-r)+n-k+2\} = (1/j^{k-2})\Pi_{r=1}^{k-2}c_{r,j-2}.$$

Hence, recalling Theorem 2, upon a little simplification,

$$Q(k,n) = 1 + \frac{n!}{(n-k+2)!(n-k+2)} \Sigma_{j=2}^{n-k+2} (-1)^{j+1}\{j\binom{n-k+2}{j}\}(\Pi_{r=1}^{k-2}c_{r,j-2})^{-1} = 1 - \overline{Q}(k,n).$$

## 4. Problem III: Random selection of $k$ pieces

Suppose $k$ of the $n$ pieces are chosen at random, the selection being equiprobable over all such choices. Let $P_{\text{rand}}(k,n)$ be the probability that the $k$ pieces so chosen do not form a $k$-gon. Interestingly, as Theorem 3 below shows, $P_{\text{rand}}(k,n)$ does not depend on $n$, but equals $P(k,k)$, i.e., the probability of not being able to form a $k$-gon when a stick is broken at random into $k$ pieces.

**Theorem 3**. $P_{\text{rand}}(k,n) = P(k,k)$.

*Proof*. We first condition on the $k$ pieces chosen and then employ Proposition 1 to get

$$P_{\text{rand}}(k,n) = \{1/\tbinom{n}{k}\}\Sigma_V P\{Y_{(i_k)} > Y_{(i_1)} + \ldots + Y_{(i_{k-1})}\}, \qquad (11)$$

where $V = \{(i_1,\ldots,i_k): i_1,\ldots,i_k$ are integers satisfying $1 \le i_1 < \ldots < i_k \le n\}$ and $\Sigma_V$ denotes sum over $(i_1,\ldots,i_k) \in V$. Recalling the well-known form of the joint density of $Y_{(i_1)},\ldots,Y_{(i_k)}$,

$$P\{Y_{(i_k)} > Y_{(i_1)} + \ldots + Y_{(i_{k-1})}\}$$

$$= \frac{n!}{\Pi_{r=1}^{k+1}(i_r - 1 - i_{r-1})!} \int_T \{\Pi_{r=1}^{k}\psi(z_r)\}[\Pi_{r=1}^{k+1}\{\Psi(z_r)-\Psi(z_{r-1})\}^{i_r-1-i_{r-1}}]dz, \qquad (12)$$

where $dz \equiv dz_1 \ldots dz_k$, $T = \{(z_1,\ldots,z_k)$: $0 < z_1 < \ldots < z_{k-1}$ and $z_k > \Sigma_{r=1}^{k-1} z_r\}$, $\psi(.)$ and $\Psi(.)$ denote, respectively, the density and distribution functions of the exponential distribution with mean 1, and

$$i_0 = 0, \quad i_{k+1} = n+1, \quad \Psi(z_0) = 0, \quad \Psi(z_{k+1}) = 1. \qquad (13)$$

Note that the right-hand side of (12) is written in terms of dummy variables $z_1,\ldots,z_k$, rather than the more customary $y_{(i_1)},\ldots,y_{(i_k)}$. This not only simplifies the notation but also makes the derivation more transparent, as we shall now see. Write $j_r = i_r - 1 - i_{r-1}$, $r = 1,\ldots, k$. Then using (13),

$$\Sigma_V \frac{(n-k)!}{\Pi_{r=1}^{k+1}(i_r-1-i_{r-1})!}\Pi_{r=1}^{k+1}\{\Psi(z_r)-\Psi(z_{r-1})\}^{i_r-1-i_{r-1}}$$

$$=\Sigma\frac{(n-k)!}{(\Pi_{r=1}^{k}j_r!)(n-k-j_1-\ldots-j_k)!}[\Pi_{r=1}^{k}\{\Psi(z_r)-\Psi(z_{r-1})\}^{j_r}]\{1-\Psi(z_k)\}^{n-k-j_1-\ldots-j_k},$$

where the sum on the right-hand side extends over integers $j_1,\ldots,j_k \geq 0$ such that $j_1+\ldots+j_k \leq n-k$, and hence equals 1. Therefore, by (11) and (12), together with a comparison with (2),

$$P_{\text{rand}}(k,n) = \{1/\tbinom{n}{k}\}\tfrac{n!}{(n-k)!}\int_T\{\Pi_{r=1}^{k}\psi(z_r)\}dz = k!\int_T \exp(-\Sigma_{r=1}^{k}z_k)dz = P(k,k),$$

which proves the result. □

One may wonder about a possible alternative proof of Theorem 3 by considering the marginal joint distribution of a normalized version of the lengths of any $k$ of the $n$ pieces. While this may have some intuitive appeal, a formalization of such a proof, with necessary arguments for all intermediate steps, is unlikely to be any simpler or shorter than the one presented above.

Let $H(k,n)$ be the number of choices of $k$ pieces out of $n$ that do not form a $k$-gon when a stick is randomly broken into $n$ pieces, $3 \leq k \leq n$. As a corollary to Theorem 3, we readily get an expression for $E\{H(k,n)\}$. To see this, note that $H(k,n) = \Sigma_V H\{(i_1,\ldots,i_k);n\}$, where $H\{(i_1,\ldots,i_k);n\}$ equals 1 if $\Delta_{(i_k)} > \Delta_{(i_1)}+\ldots+\Delta_{(i_{k-1})}$ and 0, otherwise, so that by Proposition 1,

$$E[H\{(i_1,\ldots,i_k);n\}] = P\{Y_{(i_k)} > Y_{(i_1)}+\ldots+Y_{(i_{k-1})}\}.$$

Hence, by (11) and Theorem 3,

$$E\{H(k,n)\} = \tbinom{n}{k}P_{\text{rand}}(k,n) = \tbinom{n}{k}P(k,k) = \tbinom{n}{k}k/2^{k-1},$$

because $P(k,k) = k/2^{k-1}$, as noted earlier in Section 2.

**Appendix: Reconciliation**

We show that the expression for $P(k,n)$ in Theorem 1 is in agreement with its counterpart in Verreault (2022b). For ease in presentation, this is done through a few steps.

*Proof of a recursion relation*: We first prove the recursion relation in (5) and (6). Write (4) in matrix notation as $x = Ay$, where $x = (x_1,\ldots, x_n)'$ and $y = (y_1,\ldots, y_n)'$. Clearly, then the matrix with rows $b_1',\ldots,b_n'$, as given by (1), equals $A^{-1}$. Hence, writing $1_n$ for the $n$x1 vector of ones.

$$(\beta_1,\ldots,\beta_n) = \Sigma_{r=1}^{n}b_r' = 1_n' A^{-1}, \quad \text{i.e.,} \quad (\beta_1,\ldots,\beta_n)A = 1_n'. \qquad \text{(A.1)}$$

Now, let $A = (a_{rj})$. By (4), the following hold for $A$:

(i) each diagonal element of $A$ equals 1,

(ii) any off-diagonal element $a_{rj}$ equals –1 if

either $j = r-1, 2 \leq r \leq n$ or $r-k+1 \leq j \leq r-2, k \leq r \leq n$;

i.e., interchanging the order of $j$ and $r$, if

either $r = j + 1,\ 1 \le j \le n - 1$ or $\max(k, j+2) \le r \le \min(n, j+k-1),\ 1 \le j \le n-2$,

or equivalently, splitting the range of $j$, if either $r = j + 1,\ 1 \le j \le k - 3$ or $k \le r \le \min(n, j+k-1),\ 1 \le j \le k-3$,

or $j + 1 \le r \le \min(n, j+k-1),\ k-2 \le j \le n-1$; (A.2)

else $a_{rj} = 0$.

The first two conditions in (A.2) arise only when $k \ge 4$. By (A.1), and (i), (ii) above,

$$\beta_n = 1, \quad \beta_j = 1 + \Sigma_{r=j+1}^{\min(n,j+k-1)} \beta_r \quad (j = k-2, \dots, n-1),$$

$$\beta_j = 1 + \beta_{j+1} + \Sigma_{r=k}^{\min(n,j+k-1)} \beta_r \quad (j = 1, \dots, k-3).$$

The truth of (5) is now immediate invoking (6).

*Some notation*: We now tune our notation with that in Verreault (2022b) for the purpose of the desired reconciliation. Let $\xi_j = \beta_{n+k-2-j}$, $0 \le j \le n$, and if $k \ge 4$, then let $\mu_j = \beta_{k-1-j}$, $2 \le j \le k-2$. Our result in Theorem 1 can then be expressed as

$$P(k,n) = n!/\{(\Pi_{j=k-2}^{n} \xi_j)(\Pi_{j=2}^{k-2} \mu_j)\}, \tag{A.3}$$

where the product $\Pi_{j=2}^{k-2} \mu_j$ does not appear when $k = 3$. Moreover, by (5) and (6),

$$\xi_0 = \dots = \xi_{k-3} = 0, \quad \xi_{k-2} = 1, \quad \xi_j = 1 + \Sigma_{r=1}^{k-1} \xi_{j-r} \ (j = k-1, \dots, n), \tag{A.4}$$

and if $k \ge 4$, then

$$\mu_2 = 1 + \xi_n + \Sigma_{r=2}^{k-2} \xi_{n-r}, \quad \mu_j = 1 + \mu_{j-1} + \Sigma_{r=2}^{k-j} \xi_{n-r} \ (3 \le j \le k-2), \tag{A.5}$$

with the equation for $3 \le j \le k-2$ arising in (A.5) only when $k \ge 5$.

*Verreault's result*: For completeness, we now present Verreault's (2022b) result on $P(k,n)$. Let

$$F_u = 0\ (u = 0, \dots, k-3), \quad F_{k-2} = 1, \quad F_u = \Sigma_{r=1}^{k-1} F_{u-r}\ (u \ge k-1), \tag{A.6}$$

$$f(j) = 0\ (j = 0, \dots, k-3), \quad f(j) = \Sigma_{r=k-2}^{j} F_r\ (j \ge k-2), \tag{A.7}$$

$$g(u) = 1 + \Sigma_{r=2}^{u} f(n-r)\ (u = 2, \dots, k-2), \quad h(j) = f(n) + \Sigma_{r=2}^{j} g(k-r)\ (j = 2, \dots, k-2). \tag{A.8}$$

Then Verreault (2022b) obtained $P(k,n)$ as

$$P(k,n) = n!/[\{\Pi_{j=k-2}^{n} f(j)\}\{\Pi_{j=2}^{k-2} h(j)\}], \tag{A.9}$$

where the product $\Pi_{j=2}^{k-2} h(j)$ does not appear when $k = 3$. One notes a striking similarity between (A.3) and (A.9) which will now be formalized.

*Simplifying features and reconciliation*: We now note that the $f(j)$ and $h(j)$ in (A.7) and (A.8) can be expressed in a simpler manner as

$$f(0) = \ldots = f(k-3) = 0, \quad f(k-2) = 1, \quad f(j) = 1 + \Sigma_{r=1}^{k-1} f(j-r) \ (j = k-1,\ldots, n), \qquad \text{(A.10)}$$

$$h(2) = 1 + f(n) + \Sigma_{r=2}^{k-2} f(n-r), \qquad h(j) = 1 + h(j-1) + \Sigma_{r=2}^{k-j} f(n-r) \ (3 \le j \le k-2), \quad \text{(A.11)}$$

where the $F_u$ and $g(u)$ do not appear. As in (A.5), the two equations in (A.11) are meaningful only when $k \ge 4$ and $k \ge 5$, respectively. All equations in (A.10) and (A.11), except the last one in (A.10), follow readily from (A.6)-(A.8). See Verreault (2022b, Lemma 3) for an inductive proof of the last equation in (A.10); for a direct proof, one has to note that $f(j) = \Sigma_{u=0}^{j} F_u$ ($j \ge 0$), substitute this on the right-hand side with $j$ replaced by $j - r$ and then change the order of summation.

A comparison between (A.10), (A.11) with (A.4), (A.5) yields $\xi_j = f(j)$ ($j = k - 2,\ldots, n$), and $\mu_j = h(j)$ ($2 \le j \le k - 2$), when $k \ge 4$. Hence our expression for $P(k,n)$ in (A.3) is in agreement with Verreault's (2022b) result as shown in (A.9).

**Acknowledgement**. This work was supported by a grant from the Science and Engineering Research Board, Government of India.